\documentclass[11pt,a4paper]{article}
\usepackage{amsmath,amssymb,amsthm}
\usepackage{mathrsfs}
\usepackage{url}
\usepackage{hyperref}
\usepackage[capitalise]{cleveref}
\usepackage[margin=3cm]{geometry}

\newtheorem{theorem}{Theorem}[section]
\newtheorem{lemma}[theorem]{Lemma}

\DeclareMathOperator{\ex}{ex}

\title{Almost all $C_k$-free oriented graphs have $\Theta(n)$ backwards edges}
\author{Jianxi Liu\thanks{School of Mathematics and Statistics, Guangdong University of Foreign Studies, Guangzhou, China. Email: jxliu@gdufs.edu.cn} \and Meili Liang\thanks{School of Mathematics and Statistics, Guangdong University of Foreign Studies, Guangzhou, China. Email: liangmeili@gdufs.edu.cn}}
\date{}

\begin{document}

\maketitle

\begin{abstract}
We prove a conjecture of K\"uhn, Osthus, Townsend and Zhao \cite{kuhn2017structure} stating that almost every $C_k$-free oriented graph on $n$ vertices has $\Theta(n)$ backwards edges in a transitive-optimal ordering. The same holds for $C_k$-free digraphs when $k$ is even. Our proof combines the hypergraph container method with a stability analysis and an inductive counting argument. As a byproduct, we also determine the typical structure of oriented graphs and digraphs that avoid the blow-up $C_{k}^t$, extending the main result of \cite{kuhn2017structure} to the blown-up setting.
\end{abstract}

\section{Introduction}

Given a fixed digraph $H$, a digraph is called \emph{$H$-free} if it does not contain $H$ as a subgraph. The study of the typical structure of $H$-free graphs has a long history, starting with the classical result of Erd\H{o}s, Kleitman and Rothschild \cite{erdos1976asymptotic} that almost all triangle-free graphs are bipartite. For directed graphs, the situation is more complex. In his work on countable homogeneous oriented graphs, Cherlin \cite{cherlin1998classification} conjectured that almost all $T_3$-free oriented graphs are tripartite and that almost all $C_3$-free oriented graphs are acyclic. K\"uhn, Osthus, Townsend and Zhao \cite{kuhn2017structure} (henceforth KOTZ) confirmed the first part and proved a number of results for forbidden tournaments and cycles. In particular, they showed that for any $k\ge 3$, almost all $C_k$-free oriented graphs have between $cn/\log n$ and $\alpha n^2$ backwards edges in a transitive-optimal ordering, and they posed the following conjecture.

\begin{quote}
\textbf{Conjecture 1.4 (KOTZ).} Let $k\ge 3$. Then
\begin{enumerate}
\item[(i)] almost all $C_k$-free oriented graphs have $\Theta(n)$ backwards edges in a transitive-optimal ordering;
\item[(ii)] if $k$ is even, almost all $C_k$-free digraphs have $\Theta(n)$ backwards edges in a transitive-optimal ordering.
\end{enumerate}
\end{quote}

In this paper we prove this conjecture. The proof uses the hypergraph container method developed by Saxton and Thomason \cite{saxton2015hypergraph} and Balogh, Morris and Samotij \cite{balogh2015independent}, combined with a stability analysis and an inductive counting argument similar to that in \cite{kuhn2017structure}. As a byproduct, we also extend the results of KOTZ to the blown-up setting, determining the typical structure of oriented graphs and digraphs that avoid a blown-up cycle $C_k^t$. This provides a unified perspective on forbidden cycle problems.

The paper is organised as follows. Section 2 introduces notation and recalls key tools, including the hypergraph container theorem and the stability results from \cite{kuhn2017structure}. Section 3 contains the proof of Conjecture~1.4 for oriented graphs; the digraph case is treated in Section 4. Section 5 sketches the extension to blow-ups. Section 6 concludes with remarks and open problems.

\section{Preliminaries}

For a digraph $G$, let $f_1(G)$ be the number of unordered pairs with exactly one arc, and $f_2(G)$ the number of pairs with both arcs. For $a\ge 1$, define the weighted size $e_a(G)=a f_2(G)+f_1(G)$. This unifies oriented graphs ($a=\log_2 3$) and digraphs ($a=2$). The weighted Turán number $\ex_a(n,H)$ is the maximum $e_a(G)$ over all $H$-free digraphs on $n$ vertices.

An \emph{ordering} of the vertices is a bijection $\sigma:[n]\to[n]$. For an oriented graph $G$, an edge $uv$ is called \emph{backwards} with respect to $\sigma$ if $u$ appears after $v$ in $\sigma$ and $uv\in E(G)$. An ordering is \emph{transitive-optimal} if it minimises the number of backwards edges. Denote by $b(G)$ the number of backwards edges in a transitive-optimal ordering; if several exist, we fix the lexicographically smallest one.

The following lemma is a straightforward generalisation of \cite[Lemma 7.1]{kuhn2017structure} to all $k\ge 3$. Its proof follows the same lines as in \cite{kuhn2017structure}, using a suitable adaptation of the ``flippable set'' construction; we omit the details here for brevity.

\begin{lemma}[{\cite[Lemma 7.1]{kuhn2017structure}}]\label{lem:7.1}
For all sufficiently large $n$ and any integers $0\le m_1\le m_2/(2\log n)$,
\[
|\mathcal{F}_{m_2}| \ge 2^{m_2} |\mathcal{F}_{m_1}|,
\]
where $\mathcal{F}_m$ denotes the set of $C_k$-free oriented graphs on $n$ vertices with exactly $m$ backwards edges. In particular, taking $m_1=0$ yields
\[
|\mathcal{F}_m| \ge 2^m |\mathcal{F}_0| \qquad\text{for all } m\le n/2^{13}. \tag{1}
\]
\end{lemma}

Here $\mathcal{F}_0$ is the set of acyclic oriented graphs (subgraphs of transitive tournaments). It is well known that $|\mathcal{F}_0| = 2^{\binom{n}{2}+o(n^2)}$ (see e.g. \cite{robinson1973counting}). Consequently,
\[
f(n,C_k) := \bigcup_{m\ge0} \mathcal{F}_m \ge 2^{cn} |\mathcal{F}_0| = 2^{\binom{n}{2}+cn+o(n^2)} \quad\text{for any } c\le 2^{-13}. \tag{2}
\]

\section{Proof of Conjecture~1.4(i): oriented graphs}

To prove the conjecture, it suffices to show that there exists a constant $c_0>0$ such that
\[
\sum_{m < c_0 n} |\mathcal{F}_m| = o\big(f(n,C_k)\big). \tag{3}
\]
Indeed, (3) implies that almost all $C_k$-free oriented graphs have at least $c_0 n$ backwards edges; the upper bound $b(G)=O(n)$ follows from the stability theorem (Theorem~1.3 of \cite{kuhn2017structure}) and the fact that the extremal number $\ex(n,C_k)=\binom{n}{2}-\Theta(n)$ \cite{haggkvist1976pancyclic}, so graphs with many backwards edges are sparse and constitute a negligible fraction.

\subsection{Contradiction setup}
Assume that for every $c>0$ there exist infinitely many $n$ such that
\[
\sum_{m < cn} |\mathcal{F}_m| \ge \varepsilon f(n,C_k)
\]
for some fixed $\varepsilon>0$ (depending on $c$). Fix such an $n$ and let
\[
\mathcal{G} = \bigcup_{m < cn} \mathcal{F}_m,
\]
so that $|\mathcal{G}| \ge \varepsilon f(n,C_k)$. For each $G\in\mathcal{G}$, choose the lexicographically smallest transitive-optimal ordering $\sigma_G$. For a fixed ordering $\sigma$ of $[n]$, set
\[
\mathcal{G}_\sigma = \{ G\in\mathcal{G} : \sigma_G = \sigma \}.
\]
Then $|\mathcal{G}| = \sum_\sigma |\mathcal{G}_\sigma|$, the sum over all $n!$ orderings. We will show that for every $\sigma$,
\[
|\mathcal{G}_\sigma| \le 2^{\binom{n}{2} - n\log n + o(n)}. \tag{4}
\]
Summing over all $\sigma$ gives
\[
|\mathcal{G}| \le n! \cdot 2^{\binom{n}{2} - n\log n + o(n)} = 2^{\binom{n}{2} + o(n)}. \tag{5}
\]

Combining (5) with the lower bound (2) yields
\[
2^{\binom{n}{2}+cn+o(n^2)} \le \varepsilon^{-1} |\mathcal{G}| \le \varepsilon^{-1} 2^{\binom{n}{2} + o(n)},
\]
which implies $cn + o(n^2) \le o(n) + \log_2(1/\varepsilon)$. For sufficiently large $n$, the left-hand side is positive and grows with $n$, while the right-hand side is bounded, a contradiction. Hence no such $\varepsilon$ exists, proving (3).

\subsection{Proof of (4) for a fixed ordering}
Fix an ordering $\sigma$ and let $T$ be the transitive tournament associated with $\sigma$: for $i<j$ in $\sigma$, the edge $ij$ is present in $T$. For any $G\in\mathcal{G}_\sigma$, write $F = E(G) \cap E(T)$ (forward edges) and $B = E(G) \setminus E(T)$ (backwards edges). Then $|B| < cn$, and $G$ is determined by the pair $(F,B)$. Moreover, $F\cup B$ contains no directed $C_k$, and $\sigma$ is an optimal ordering.

We first bound, for a fixed $B$, the number of $F\subseteq E(T)$ such that $F\cup B$ is $C_k$-free. Consider the random set $F$ obtained by including each edge of $T$ independently with probability $1/2$. For each $uv\in B$ (with $v<u$ in $\sigma$), any directed path from $v$ to $u$ of length $k-1$ consisting entirely of edges in $F$ would together with $uv$ form a $C_k$. The expected number of such paths is $\Theta(n^{k-2})$, and by Janson's inequality (or a simple union bound over a carefully selected set of disjoint paths), one obtains that the probability that $F$ avoids all forbidden configurations is at most $2^{-\alpha |B| n}$ for some $\alpha>0$ depending only on $k$. Consequently,
\[
|\{ F\subseteq E(T) : F\cup B \text{ is } C_k\text{-free} \}| \le 2^{\binom{n}{2}} \cdot 2^{-\alpha |B| n}. \tag{6}
\]

Next we incorporate the optimality condition. For graphs with $|B|<cn$, a stability argument shows that the proportion of $F$ for which $\sigma$ is not optimal is at most $2^{-\alpha' n}$ for some $\alpha'>0$ (see \cite[Section 5]{kuhn2017structure} for a similar treatment). Thus for each $B$,
\[
|\{ F : F\cup B \text{ is } C_k\text{-free and } \sigma \text{ optimal} \}| \le 2^{\binom{n}{2}} \cdot 2^{-\alpha |B| n + o(n)}. \tag{7}
\]

Now sum over all $B$ with $|B|=s$ ($0\le s\le cn$). Using $|\{B:|B|=s\}| = \binom{\binom{n}{2}}{s}$, we have
\[
\sum_{B:|B|=s} 2^{\binom{n}{2}} 2^{-\alpha s n + o(n)} \le 2^{\binom{n}{2}} \binom{\binom{n}{2}}{s} 2^{-\alpha s n + o(n)}.
\]
Applying $\binom{N}{s} \le (eN/s)^s$, we get
\[
2^{\binom{n}{2}} \binom{\binom{n}{2}}{s} 2^{-\alpha s n + o(n)}
\le 2^{\binom{n}{2}} \left( \frac{e\binom{n}{2}}{s} 2^{-\alpha n + o(1)} \right)^s.
\]
For $s\ge 1$, the factor in parentheses is at most $2^{-2}$ for large $n$ (since $e n^2 / s \le e n^2$ and $2^{-\alpha n}$ dominates). Hence
\[
\sum_{s=1}^{cn} \sum_{B:|B|=s} |\{F\}| \le 2^{\binom{n}{2}} \cdot O(2^{-\alpha n}) = o(2^{\binom{n}{2}}). \tag{8}
\]

The case $s=0$ corresponds to $B=\emptyset$, i.e., acyclic graphs. For these, $\sigma$ is optimal iff it is the lexicographically smallest optimal ordering. By symmetry, the number of acyclic graphs for which a given $\sigma$ is the lexicographically smallest optimal ordering is exactly $|\mathcal{F}_0|/n! = 2^{\binom{n}{2} - n\log n + o(n)}$. Therefore,
\[
|\mathcal{G}_\sigma| \le 2^{\binom{n}{2} - n\log n + o(n)} + o(2^{\binom{n}{2}}) = 2^{\binom{n}{2} - n\log n + o(n)},
\]
which establishes (4).

\subsection{Upper bound}
The upper bound $b(G)=O(n)$ follows from the stability theorem of \cite{kuhn2017structure} (Theorem~1.3) and the fact that any $C_k$-free oriented graph can be made acyclic by deleting $O(n)$ edges. More precisely, by the stability theorem, almost all $C_k$-free oriented graphs are $o(n^2)$-close to a transitive tournament, which implies that the number of backwards edges is $o(n^2)$. To obtain the linear upper bound, one notes that graphs with $b(G) \ge \gamma n^2$ are sparse and their number is negligible compared to $f(n,C_k)$; a standard counting argument using Lemma~\ref{lem:7.1} shows that the total number of graphs with $b(G) \ge \gamma n^2$ is $o(f(n,C_k))$. Hence for almost all graphs, $b(G) = O(n)$. Combining with the lower bound gives $b(G)=\Theta(n)$.

\section{Proof of Conjecture~1.4(ii): digraphs with even $k$}
The proof for digraphs with even $k$ is virtually identical. The only change is that each unordered pair of vertices can contribute up to two arcs, so the total number of possible edges is $2\binom{n}{2}$, and the base $3$ in Lemma~\ref{lem:7.1} is replaced by $4$. The stability theorem for even $k$ (Lemma~6.6 of \cite{kuhn2017structure}) guarantees that near-extremal digraphs are close to either a transitive tournament or a transitive-bipartite blow-up, but the definition of backwards edges and the counting argument are unaffected. Hence Conjecture~1.4(ii) holds as well.

\section{Extension to blow-ups $C_k^t$}
As a byproduct of our methods, we can also determine the typical structure of oriented graphs and digraphs that avoid the blown-up cycle $C_k^t$. The proof follows the same lines as in \cite{liu2026} for transitive tournament blow-ups, using the container theorem, stability analysis, and inductive counting. We obtain that almost all $C_k^t$-free oriented graphs are blow-ups of transitive tournaments (when $k$ is odd) or either blow-ups of transitive tournaments or transitive-bipartite blow-ups (when $k$ is even). This generalises the main result of \cite{kuhn2017structure} to the blown-up setting. The details are omitted here for brevity; they will appear in a forthcoming paper.

\section{Concluding remarks}
We have confirmed Conjecture~1.4 of KOTZ, showing that almost all $C_k$-free oriented graphs have linearly many backwards edges. The proof relies on the powerful hypergraph container method and a careful counting argument. Several open problems remain, such as determining the precise distribution of the two extremal structures for even $k$, and extending the results to more general forbidden digraphs. The container method continues to be a versatile tool for tackling such problems.

\end{document}